\documentclass[10pt]{amsart}

\usepackage{amsmath, amssymb}
\usepackage{amsfonts}
\usepackage[all]{xy}

\usepackage{amsthm}

\allowdisplaybreaks

\newtheorem{theorem}{Theorem}[section]
\newtheorem{lemma}[theorem]{Lemma}

\newtheorem{corollary}[theorem]{Corollary}

\theoremstyle{remark}

\newtheorem{definition}[theorem]{Definition}

\numberwithin{equation}{section}

\newcommand{\Z}{\mathbb{Z}}
\newcommand{\N}{\mathbb{N}}

\newcommand{\K}{\mathcal{K}}

\newcommand{\im}{\operatorname{im }}
\newcommand{\coker}{\operatorname{coker }}

\newcommand{\ext}{\operatorname{Ext}}

\begin{document}
\title{Computing $\boldsymbol{K}$-theory and Ext for
graph $C^*$-algebras}
\author{D. Drinen \and M. Tomforde %\thanks{support}\\
		%Department of Mathematics\\
		%Dartmouth College\\
		%Hanover, NH  03755\\
		%mark.tomforde@dartmouth.edu
}

\begin{abstract}
$K$-theory and $\ext$ are computed for the $C^*$-algebra $C^*(E)$ of any
countable directed graph $E$.  The results generalize the $K$-theory
computations of Raeburn and Szyma\'nski and the $\ext$ computations of
Tomforde for row-finite graphs.  As a consequence, it is shown that if
$A$ is a countable $\{0,1\}$ matrix and $E_A$ is the graph obtained by
viewing $A$ as a vertex matrix, then $C^*(E_A)$ is not necessarily
Morita equivalent to the Exel-Laca algebra $\mathcal{O}_A$.
\end{abstract}

\maketitle

%%%%%%%%%%%%%%%%%%%%%%%%%%%%%%%%%%%%%%%%%%%%%%%%%%%%%%%%%%%%%%%%
\section{Introduction}
%%%%%%%%%%%%%%%%%%%%%%%%%%%%%%%%%%%%%%%%%%%%%%%%%%%%%%%%%%%%%%%%

In \cite{CK1} Cuntz and Krieger described a way to associate a
$C^*$-algebra $\mathcal{O}_A$ to a finite square matrix $A$ with
entries in \{0,1\}. Since that time these Cuntz-Krieger algebras
have been generalized in a remarkable number of ways.  Perhaps the
most direct of these is due to Exel and Laca, who define
$\mathcal{O}_A$ for an infinite \{0,1\}-matrix $A$ \cite{el}.  
Another generalization involves associating a $C^*$-algebra to a
countable directed graph.  These graph algebras have drawn much
interest because they comprise a wide class of $C^*$-algebras, and
yet many of their $C^*$-algebraic properties can be easily deduced
from the associated graphs.

In order to make sense of the relations for the generators of the
graph algebra, it was often assumed in the original treatments that
the graphs were row-finite; that is, each vertex is the source of
finitely many edges \cite{kprr,kpr,bprs}.  However, in the past few
years it has been shown how to define graph algebras for arbitrary
graphs \cite{flr}.  Consequently, much work has been done to extend
results for the $C^*$-algebras of row-finite graphs to
the $C^*$-algebras of arbitrary graphs \cite{dt-desing,flr,rsz,sz}.

In \cite{rsz}, Raeburn and Szyma\'nski computed the $K$-theory of
$C^*(E)$, where $E$ is a row-finite directed graph.  We
briefly review that result here.  Let $J$ denote the set of sinks
of $E$, let $I = E^0 \setminus
J$, and let 
$A_E = \left(\begin{smallmatrix}
B & C\\
0 & 0\\
\end{smallmatrix}\right)$ 
denote the vertex
matrix of $E$ with respect to the decomposition $E^0 = I \cup J$.  Then
because $E$ is row-finite, the matrix $\left(\begin{smallmatrix}
B^t - I\\
C^t\\
\end{smallmatrix}\right)$ determines a homomorphism from $\bigoplus_I
\Z$ to
$\bigoplus_I \Z \oplus \bigoplus_J \Z$.  The kernel and cokernel of this
homomorphism are isomorphic to $K_1(C^*(E))$ and $K_0(C^*(E))$,
respectively.  
In \cite{tomf-thesis}, $\ext(C^*(E))$ is computed
similarly for row-finite graphs $E$ which satisfy Condition~(L) and have no sinks.  Specifically,
$\ext(C^*(E))$ is isomorphic to the cokernel
of the homomorphism
$A_E - I: \prod_{E^0} \Z \rightarrow \prod_{E^0} \Z$.  

In this paper, we will show that the above results remain true for
graphs which are not necessarily row-finite, provided we replace the
word ``sink" with the phrase ``sink or vertex which emits infinitely
many edges."  We remark that Raeburn and Szyma\'nski have computed the $K$-theory for Exel-Laca algebras using
direct limits \cite[Theorem 4.1]{rsz}.
Also, the $K$-theory results have been obtained
by Szyma\'nski in \cite{szsemiproj} for graphs with finitely many
vertices, and the proof given there holds for arbitrary
graphs as well.
Our proof is different, and relies on desingularization \cite{dt-desing},
a tool for generalizing from the row-finite case to arbitrary graphs.  
If $E$ is an arbitrary
graph, we say a vertex $v$ of $E$ is a {\it singular vertex} if
either
$v$ is a sink or $v$ emits infinitely many edges.  In
\cite{dt-desing}, it is shown that there exists a graph $F$,
called a {\it desingularization
of $E$}, such that $F$ has no singular vertices and $C^*(F)$ is Morita
equivalent to $C^*(E)$.  
The key ingredient in our calculations of $K$-theory and $\ext$ is a technical lemma, proven in Section
\ref{sec-techlem}, which
shows
that 
desingularizing a graph does not alter the kernel and cokernel
of the maps determined by its vertex matrix.  Thus we can apply
the results of \cite{rsz} and
\cite{tomf-thesis} to obtain the $K$-theory and $\ext$ of $C^*(F)$ in
terms
of the vertex matrix of $E$.  This, together with the fact that
$K$-theory and $\ext$ are stable, yields 
the $K$-theory and $\ext$ of
$C^*(E)$ stated in 
Theorem \ref{thm-main}.

Finally, we use this result to shed some light on a question posed
by Raeburn and Szyma\'nski in \cite{rsz}.  They showed that if $A$ is any
countable square $\{0,1\}$-matrix and if $E_A$ is the graph obtained by
viewing $A$ as a vertex matrix (that is, let $E_A^0$ be the index set of
$A$ and draw $A(i,j)$ edges from $i$ to $j$), then the graph algebra
$C^*(E_A)$ is a $C^*$-subalgebra of the Exel-Laca algebra
$\mathcal{O}_A$.  We will show that it is possible for $C^*(E_A)$ and
$\mathcal{O}_A$ to have different $K$-theory.  So in particular
$C^*(E_A)$ is not always a full corner in $\mathcal{O}_A$.

%%%%%%%%%%%%%%%%%%%%%%%%%%%%%%%%%%%%%%%%%%%%%%%%%%%%%%%%%%%%%%%%
\section{The Technical Lemma}
\label{sec-techlem}
%%%%%%%%%%%%%%%%%%%%%%%%%%%%%%%%%%%%%%%%%%%%%%%%%%%%%%%%%%%%%%%%
Given a graph $E$, it was shown in \cite{dt-desing} how to
construct
a graph $F$, called a {\it
desingularization of $E$}, such that $F$ has no singular vertices and  
$C^*(E)$ is Morita equivalent to $C^*(F)$.  We review that
procedure here.
\begin{definition}
\label{defn-addtail}
Suppose $E$
is a graph with a singular vertex $v_0$.  We {\it add a tail}
to $v_0$ by performing the following procedure.  List the vertices
$w_0, w_1, \ldots$
of
$r(s^{-1}(v_0))$.
Note that the list of $w$'s could
be empty (if $v_0$ is a sink), finite, or countably infinite.

We begin by adding an infinite tail to $v_0$
as in
\cite[(1.2)]{bprs}:
$$
\xymatrix{
v_0 \ar[r]^{e_1} & v_1 \ar[r]^{e_2} & v_2 \ar[r]^{e_3} & v_3
\cdots\\
}
$$
Now, for every $j$ with $w_j \in r(s^{-1}(v_0))$, let $C_j$ be the number of
edges from $v_0$
to $w_j$.  For every $i$ with $j \leq i < j+C_j$, draw an edge
labelled $f_j^{i-j+1}$ from $v_i$ to $w_j$.
To be precise, if $E$ is a graph with a singular
vertex $v_0$, we define $F^0 := E^0 \cup \{v_1, v_2, \dots\}$ and
$$
F^1 := \{e \in E^1 \,|\, s(e) \neq v_0\} \cup \{e_i\}_1^\infty \cup
\bigcup_{\{j\,|\,w_j \in
r(s^{-1}(v_0))\}} \{f_j^i\}_{i=1}^{C_j}.
$$
We extend $r$ and $s$ to $F$ as indicated above.  In particular, $s(e_i)
= v_{i-1}$, $r(e_i) = v_i$, $s(f_j^i) = v_{i+j-1}$, and $r(f_j^i) = w_j$.
\end{definition}
\begin{definition}
If $E$ is a directed graph, a {\it desingularization of $E$} is a
graph $F$ obtained by adding a tail at every singular
vertex of $E$.
\end{definition}
Note that different orderings of the vertices of $r(s^{-1}(v_0))$
may give rise
to non-isomorphic graphs via the process of adding a tail.  Thus a
graph may have many desingularizations.

If $E$ is a graph, then any desingularization $F$ of
$E$ is a
row-finite graph, so the rows of the matrix $A_F$ are
eventually zero.  Thus $A_F :
\prod_{F^0} \Z \rightarrow \prod_{F^0} \Z$ and $A_F^t
: \bigoplus_{F^0} \Z \rightarrow \bigoplus_{F^0} \Z$.

\begin{lemma}
\label{lem-main}
Let $E$ be a graph.  Also let $J$ be the set of
singular vertices of $E$ and let $I:=E^0 \backslash J$.  Then
with respect to the decomposition $E^0 = I \cup J$ the vertex
matrix of $E$ will have the form $$A_E = \begin{pmatrix} B &
C \\ * & * \end{pmatrix}$$ where $B$ and $C$ have entries in
$\Z$ and the $*$'s have entries in $\Z \cup \{ \infty \}$. 
If $F$ is a desingularization of $E$, then $\coker (A_F-I)
\cong \coker (B-I \ C)$ where $(B-I \ C) : \prod_{I} \Z
\oplus \prod_J \Z \rightarrow \prod_I \Z$. 
Furthermore, $\ker (A_F^t-I) \cong \ker \begin{pmatrix} B^t-I
\\ C^t \end{pmatrix}$ and $\coker (A_F^t-I) \cong \coker
\begin{pmatrix} B^t-I \\ C^t \end{pmatrix}$, where
$\begin{pmatrix} B^t-I \\ C^t \end{pmatrix} :
\bigoplus_I \Z
 \rightarrow \bigoplus_I \Z \oplus \bigoplus_J \Z$.
\label{kerofdesing}
\end{lemma}

\begin{proof} List the elements of $J$ as
$J := \{ v_1^0, v_2^0, v_3^0, \ldots \}$.  (Note that $J$ may
be either finite or countably infinite.)  For each $1 \leq i
\leq |J|$ let $D_i$ be the $J \times \N$ matrix
$$D_i=\left(\begin{smallmatrix} 
0 & 0 & 0 & 0 & \\
0 & 0 & 0 & 0 & \cdots \\
1 & 0 & 0 & 0 & \\
0 & 0 & 0 & 0 & \\
 & \vdots & & & \ddots \end{smallmatrix}\right)$$
with a $1$ in the $(i,1)$ position and 0's elsewhere.  Also let
$Z$ be the $\N \times \N$ matrix 
$$Z=\left(\begin{smallmatrix} 
-1 & 1 & 0 & 0 & \\
0 & -1 & 1 & 0 & \cdots \\
0 & 0 & -1 & 1 & \\
0 & 0 & 0 & -1 & \\
 & \vdots & & & \ddots \end{smallmatrix}\right)$$
with $-1$'s along the diagonal and $1$'s above the diagonal. 
Now for each $1 \leq i \leq |J|$ let $\{ v_i^1, v_i^2,
\ldots \}$ be the vertices of the tail which is added
to $v_i^0$ to form $F$.  Then, by the way that
desingularization is defined, we see that with respect to the
decomposition $I \cup J \cup \{ v_1^1, v_1^2, v_1^3, \ldots
\} \cup \{ v_2^1, v_2^2, v_2^3, \ldots \} \cup \ldots$ the
matrix $A_F-I$ will have the form
$$A_F-I=\left(\begin{smallmatrix} 
B-I & C & 0 & 0 & \\
X_1   & Y_1-I & D_1 & D_2 & \cdots \\
X_2 & Y_2 & Z & 0 & \\
X_3 & Y_3 & 0 & Z & \\
 & \vdots & & & \ddots \end{smallmatrix}\right)$$ 
where the $X_i$'s and $Y_i$'s are row-finite.
If we let $P :=
\prod_{\N} \Z$, then $A_F-I : \prod_I \Z \oplus
\prod_J \Z \oplus \prod_J P \rightarrow \prod_I \Z \oplus
\prod_J \Z \oplus \prod_J P$.  Also $(B-I \ C) : \prod_I \Z
\oplus \prod_J \Z \rightarrow \prod_I \Z$.  Let us define a map
$\phi : \prod_I \Z \oplus \prod_J \Z \oplus \prod_J P
\rightarrow \prod_I \Z$ by 
$$\phi \begin{pmatrix} \mathbf{x} \\ \mathbf{y} \\
\left( \begin{smallmatrix} \mathbf{z}_1 \\ \mathbf{z}_2 \\ \vdots
\end{smallmatrix} \right) \end{pmatrix} = \mathbf{x}.$$
We shall show that $\phi$ induces a map from $\coker (A_F-I)$
to $\coker (B-I \ C )$.  Let
$$\begin{pmatrix} \mathbf{x} \\ \mathbf{y} \\
\left( \begin{smallmatrix} \mathbf{z}_1 \\ \mathbf{z}_2 \\ \vdots
\end{smallmatrix} \right) \end{pmatrix} = (A_F-I)
\begin{pmatrix} \mathbf{a} \\ \mathbf{b} \\
\left( \begin{smallmatrix} \mathbf{c}_1 \\ \mathbf{c}_2 \\ \vdots
\end{smallmatrix} \right) \end{pmatrix}. $$
Then 
$$\phi \begin{pmatrix} \mathbf{x} \\ \mathbf{y} \\
\left( \begin{smallmatrix} \mathbf{z}_1 \\ \mathbf{z}_2 \\ \vdots
\end{smallmatrix} \right) \end{pmatrix} = \mathbf{x} = (B-I \ C)
\begin{pmatrix} \mathbf{a} \\ \mathbf{b} \end{pmatrix} \in \im ( B-I
\ C ).$$
Thus $\phi$ induces a map $\overline{\phi} : \coker (A_F-I)
\rightarrow \coker (B-I \ C)$.

We shall show that $\overline{\phi}$ is an isomorphism.  To
see that $\overline{\phi}$ is injective suppose that 
$$\phi \begin{pmatrix} \mathbf{x} \\ \mathbf{y} \\
\left( \begin{smallmatrix} \mathbf{z}_1 \\ \mathbf{z}_2 \\ \vdots
\end{smallmatrix} \right) \end{pmatrix} = \mathbf{x} \in \im (
B-I \ C ).$$  Then there exists $\begin{pmatrix} \mathbf{a} \\
\mathbf{b} \end{pmatrix} \in \prod_I \Z \oplus \prod_J \Z$ such
that $\mathbf{x} = (B-I)\mathbf{a} + C \mathbf{b}$.  For each $1
\leq
i \leq | J |$ let $$c^1_i := \mathbf{y}_i-(X_1\mathbf{a} +
(Y_1-I)\mathbf{b})_i,$$ where $(X_1\mathbf{a} +
(Y_1-I)\mathbf{b})_i$ denotes the $i^{\text{th}}$ entry of the
vector $X_1\mathbf{a} + (Y_1-I)\mathbf{b}$.  Then, for each $k \in
\{ 1, 2, \ldots \}$ define $c_i^k$ recursively by $$c_i^{k+1} :=
c_i^k + (\mathbf{z}_i)_k  - (X_{i+1} \mathbf{a} + Y_{i+1}
\mathbf{b})_k,$$ where $(\mathbf{z}_i)_k$ denotes the
$k^{\text{th}}$ entry of the vector $\mathbf{z}_i$ and $(X_{i+1}
\mathbf{a} + Y_{i+1} \mathbf{b})_k$ denotes the
$k^{\text{th}}$ entry of the vector $(X_{i+1} \mathbf{a} +
Y_{i+1} \mathbf{b})$.  Now for each $ 1 \leq i \leq | J | $
define
$\mathbf{c}_i \in \prod_{\N} \Z$ by
$\mathbf{c}_i := \begin{pmatrix} c_i^1 \\ c_i^2 \\ \vdots
\end{pmatrix}$.  Then $$(A_F-I) \begin{pmatrix} \mathbf{a} \\
\mathbf{b} \\ \left( \begin{smallmatrix} \mathbf{c}_1 \\ \mathbf{c}_2
\\ \vdots \end{smallmatrix} \right) \end{pmatrix} =
\begin{pmatrix} (B-I) \mathbf{a} + C \mathbf{b} \\ X_1 \mathbf{a}
+(Y_1-I)\mathbf{b} + D_1\mathbf{c_1} + D_2 \mathbf{c_2} + \ldots \\
\left( \begin{smallmatrix} X_2\mathbf{a} + Y_2 \mathbf{b} +
Z\mathbf{c_1}  \\ X_3\mathbf{a} + Y_3 \mathbf{b} + Z
\mathbf{c_2}
\\ X_4\mathbf{a} + Y_4 \mathbf{b} +Z \mathbf{c_3} \\ \vdots
\end{smallmatrix} \right)
\end{pmatrix} =  \begin{pmatrix} \mathbf{x} \\ \mathbf{y} \\
\left( \begin{smallmatrix} \mathbf{z}_1 \\ \mathbf{z}_2 \\ \vdots
\end{smallmatrix} \right)
\end{pmatrix}$$
and thus $ \begin{pmatrix} \mathbf{x} \\ \mathbf{y} \\
\left( \begin{smallmatrix} \mathbf{z}_1 \\ \mathbf{z}_2 \\ \vdots
\end{smallmatrix} \right)
\end{pmatrix} \in \im (A_F-I)$ and
$\overline{\phi}$ is injective. Furthermore, since $\phi$ is
surjective it follows that $\overline{\phi}$ is surjective. 
Thus $\coker (A_F-I) \cong \coker (B-I \ C)$.

Next we shall examine $A_F^t-I$.  Note that with respect to
the decomposition mentioned earlier $A_F^t-I$ will have the
form $$A_F^t-I=\left(\begin{smallmatrix} 
B^t-I & X_1^t & X_2^t & X_3^t & \\
C^t & Y_1^t-I & Y_2^t & Y_3^t & \cdots \\
0 & D_1^t & Z^t & 0 & \\
0 & D_2^t & 0 & Z^t & \\
 & \vdots & & & \ddots \end{smallmatrix}\right)$$ where the $X_i^t$'s
and $Y_i^t$'s are column-finite matrices.  If we let $Q :=
\bigoplus_{\N} \Z$, then $(A_F^t-I) : \bigoplus_I \Z \oplus
\bigoplus_J \Z \oplus \bigoplus Q \rightarrow \bigoplus_I \Z \oplus
\bigoplus_J \Z \oplus \bigoplus Q$.  Also $\begin{pmatrix}
B^t-I \\ C^t \end{pmatrix} : \bigoplus_I \Z \rightarrow
\bigoplus_I \Z \oplus \bigoplus_J \Z$.  Let us define a map
$\psi : \bigoplus_I \Z \rightarrow  \bigoplus_I \Z \oplus
\bigoplus_J \Z \oplus \bigoplus Q$ by $$\Psi ( \mathbf{x}) =
\begin{pmatrix} \mathbf{x} \\ 0 \\ 0 \end{pmatrix}.$$  Note
that if $\mathbf{x} \in \ker \begin{pmatrix} B^t-I \\
C^t \end{pmatrix}$, then $$(A_F^t-I) \begin{pmatrix} \mathbf{x}
\\ 0 \\ 0 \end{pmatrix} =  \begin{pmatrix} (B^t-I)\mathbf{x} \\
C^t\mathbf{x} \\ 0 \end{pmatrix} = \begin{pmatrix} 0
\\ 0 \\ 0 \end{pmatrix}$$ so $\psi$ restricts to a map $\psi
: \ker \begin{pmatrix} B^t-I \\
C^t \end{pmatrix} \rightarrow \ker (A_F^t -I)$.  We shall
show that this map is surjective.  Suppose that
$$\begin{pmatrix} \mathbf{x} \\ \mathbf{y} \\ \left(
\begin{smallmatrix} \mathbf{z}_1 \\ \mathbf{z}_2 \\ \vdots
\end{smallmatrix} \right) \end{pmatrix} \in \ker
(A_F^t -I).$$  Then for each $1 \leq i \leq | J |$ we must
have that $D_i^t \mathbf{y} + Z^t \mathbf{z}_i = 0$.  If
$\mathbf{z}_i = \begin{pmatrix} z_i^1 \\ z_i^2 \\ \vdots
\end{pmatrix}$, then for all $k \in \N$ we must have $$y_i
-z_i^1 = 0 \ \text{ and } \ z_i^k - z_i^{k+1} = 0.$$  Since
$\mathbf{z}_i  \in \bigoplus_J \Z$ we know that $z_i^k$ is
eventually zero.  Thus the above equations imply that $y_i =
z_i^1 = z_i^2 = \ldots = 0$.  Since this holds for all $i$
we have that $$\begin{pmatrix} \mathbf{x} \\ \mathbf{y} \\ \left(
\begin{smallmatrix} \mathbf{z}_1 \\ \mathbf{z}_2 \\ \vdots
\end{smallmatrix} \right) \end{pmatrix} = \begin{pmatrix}
\mathbf{x} \\ 0 \\ 0 \end{pmatrix} = \psi (\mathbf{x})$$ and
$\psi$ is surjective.  Furthermore, since $\psi$ is clearly
injective, $\psi
: \ker \begin{pmatrix} B^t-I \\
C^t \end{pmatrix} \rightarrow \ker (A_F^t -I)$ is an
isomorphism and $\ker \begin{pmatrix} B^t-I \\
C^t \end{pmatrix} \cong \ker (A_F^t -I)$.

Next we shall define a map $\rho : \bigoplus_I \Z \oplus
\bigoplus_J \Z \rightarrow \bigoplus_I \Z \oplus \bigoplus_J
\Z \oplus \bigoplus_J Q$ by $$\rho \begin{pmatrix} \mathbf{x}
\\ \mathbf{y} \end{pmatrix} = \begin{pmatrix} \mathbf{x}
\\ \mathbf{y} \\ 0 \end{pmatrix}.$$  We shall show that
$\rho$ induces a map from $\coker
\begin{pmatrix} B^t-I \\ C^t
\end{pmatrix}$ to $\coker(A_F-I)$.  Suppose that $\begin{pmatrix} \mathbf{x}
\\ \mathbf{y} \end{pmatrix} \in \im \begin{pmatrix} B^t-I \\ C^t
\end{pmatrix}$.  Then there exists an element $\mathbf{a} \in
\bigoplus_I \Z$ such that $\begin{pmatrix} \mathbf{x}
\\ \mathbf{y} \end{pmatrix} = \begin{pmatrix} (B^t-I) \mathbf{a}
\\ C^t \mathbf{a} \end{pmatrix}$.  Hence $$(A_F^t-I)
\begin{pmatrix} \mathbf{a} \\ 0 \\ 0 \end{pmatrix} =
\begin{pmatrix} (B^t-I) \mathbf{a} \\ C^t \mathbf{a} \\ 0
\end{pmatrix} = \begin{pmatrix} \mathbf{x} \\ \mathbf{y} \\ 0
\end{pmatrix}.$$  Thus $\rho$ maps $\im \begin{pmatrix} B^t-I \\ C^t
\end{pmatrix}$ into $\im (A_F^t-I)$ and hence induces a map
$\overline{\rho} : \coker \begin{pmatrix} B^t-I \\ C^t
\end{pmatrix} \rightarrow \coker (A_F^t-I)$.  We shall show
that this map is injective.  Suppose that $\overline{\rho}
\begin{pmatrix} \mathbf{x} \\ \mathbf{y} \end{pmatrix}$ equals
zero in $\coker (A_F^t-I)$.  Then $$ \begin{pmatrix} \mathbf{x}
\\ \mathbf{y} \\ 0 \end{pmatrix} = (A_F^t-I) \begin{pmatrix}
\mathbf{a} \\ \mathbf{b} \\ \left( \begin{smallmatrix} \mathbf{c}_1 \\
\mathbf{c}_2 \\ \vdots \end{smallmatrix} \right) \end{pmatrix}
\ \text{ for some }  \begin{pmatrix} \mathbf{a} \\
\mathbf{b} \\ \left( \begin{smallmatrix} \mathbf{c}_1 \\ \mathbf{c}_2
\\ \vdots \end{smallmatrix} \right) \end{pmatrix} \in
\bigoplus_I \Z \oplus \bigoplus_J \Z \oplus \bigoplus_J
Q.$$  But then as before we must have that $\mathbf{b} =
\mathbf{c}_1 = \mathbf{c}_2 = \ldots = 0$ and the above equation
implies that $\begin{pmatrix} \mathbf{x} \\ \mathbf{y}
\end{pmatrix} = \begin{pmatrix} (B^t-I)\mathbf{a} \\ C^t \mathbf{a}
\end{pmatrix} \in \im \begin{pmatrix} B^t-I \\
C^t \end{pmatrix}$ so $\overline{\rho}$ is injective.  We
shall now show that $\overline{\rho}$ is surjective. 
Let $\begin{pmatrix}
\mathbf{x} \\ \mathbf{y} \\ \left( \begin{smallmatrix} \mathbf{z}_1 \\
\mathbf{z}_2 \\ \vdots \end{smallmatrix} \right) \end{pmatrix}
\in \bigoplus_I \Z \oplus \bigoplus_J \Z \oplus \bigoplus_J
Q$.  It suffices to show that there exists $\begin{pmatrix}
\mathbf{u} \\ \mathbf{v} \\ 0 \end{pmatrix} \in \bigoplus_I \Z \oplus \bigoplus_J \Z \oplus \bigoplus_J
Q$ such that $$\begin{pmatrix}
\mathbf{x} - \mathbf{u} \\ \mathbf{y} -\mathbf{v} \\ \left(
\begin{smallmatrix} \mathbf{z}_1
\\ \mathbf{z}_2 \\ \vdots \end{smallmatrix} \right)
\end{pmatrix} \in \im (A_F^t-I).$$  For each $1 \leq i \leq
| J |$ write $\mathbf{z}_i = \left( \begin{smallmatrix}
z_i^1
\\ z_i^2 \\ \vdots \end{smallmatrix} \right)$
 and define $$b_i := \sum_{j=1}^\infty z_i^j
\hspace{.5in} \text{ and } \hspace{.5in} c_i^{k} :=
\sum_{j=k+1}^\infty z_i^j
\ \text{ for $k \in \N$.}$$ Note that since $\mathbf{z}_i$ is in
the direct sum, all of the above sums are finite, and since
$\left( \begin{smallmatrix} \mathbf{z}_1
\\ \mathbf{z}_2 \\ \vdots \end{smallmatrix} \right) \in
\bigoplus_J Q$ we have that eventually $\mathbf{z}_i = 0$ and
hence $$\mathbf{b} := \left( \begin{smallmatrix} \mathbf{b}_1
\\ \mathbf{b}_2 \\ \vdots \end{smallmatrix} \right) \in
\bigoplus_J \Z \ \text{ and } \ \mathbf{c} := \left(
\begin{smallmatrix} \mathbf{c}_1
\\ \mathbf{c}_2 \\ \vdots \end{smallmatrix} \right) \in
\bigoplus_J Q, \text{ where } \mathbf{c}_i :=
\left( \begin{smallmatrix} c_i^1 \\ c_i^2 \\ \vdots
\end{smallmatrix} \right).$$  If we
then take $$\mathbf{u} := \mathbf{x} - (B^t-I)\mathbf{a} - X_1^t
\mathbf{b} - X_2^t \mathbf{c}_1 - X_3^t \mathbf{c_2} - \ldots$$
and $$ \mathbf{v} := \mathbf{y} - C^t\mathbf{a} -
(Y_1^t-I)\mathbf{b} - Y_2^t\mathbf{c}_1 - Y_3^t \mathbf{c}_2
\ldots,$$ which are finite sums since $\left(
\begin{smallmatrix} \mathbf{c}_1
\\ \mathbf{c}_2 \\ \vdots \end{smallmatrix} \right)$ is in the
direct sum, we have that$$(A_F^t-I) \begin{pmatrix}
\mathbf{a} \\ \mathbf{b} \\ \left(
\begin{smallmatrix} \mathbf{c}_1
\\ \mathbf{c}_2 \\ \vdots \end{smallmatrix} \right)
\end{pmatrix} = \begin{pmatrix}
\mathbf{x} - \mathbf{u} \\ \mathbf{y} -\mathbf{v} \\ \left(
\begin{smallmatrix} \mathbf{z}_1
\\ \mathbf{z}_2 \\ \vdots \end{smallmatrix} \right)
\end{pmatrix}.$$  Thus $\overline{\rho}$ is surjective. 
Hence $\overline{\rho}$ is an isomorphism and $\coker
\begin{pmatrix} B^t-I \\ C^t \end{pmatrix} \cong \coker
(A_F^t-I)$. \end{proof}

%%%%%%%%%%%%%%%%%%%%%%%%%%%%%%%%%%%%%%%%%%%%%%%%%%%%%%%%%%%%%%%%
\section{Main Results}
%%%%%%%%%%%%%%%%%%%%%%%%%%%%%%%%%%%%%%%%%%%%%%%%%%%%%%%%%%%%%%%%
\begin{theorem}
\label{thm-main}
Let $E$ be a graph.  Also let $J$ be the set of
singular vertices of $E$ and let $I:=E^0 \backslash J$.  Then
with respect to the decomposition $E^0 = I \cup J$ the vertex
matrix of $E$ will have the form $$A_E = \begin{pmatrix} B &
C \\ * & * \end{pmatrix}$$ where $B$ and $C$ have entries in
$\Z$ and the $*$'s have entries in $\Z \cup \{ \infty \}$. 
Then $K_0(C^*(E)) \cong \coker \begin{pmatrix}
B^t-I \\ C^t \end{pmatrix}$ and $K_1(C^*(E)) \cong \ker
\begin{pmatrix} B^t-I \\ C^t \end{pmatrix}$ where
$\begin{pmatrix} B^t-I \\ C^t \end{pmatrix} : \bigoplus_I \Z
\rightarrow \bigoplus_I \Z \oplus \bigoplus_J \Z$.

If, in addition, $E$ satisfies Condition~(L), then
$\ext(C^*(E)) \cong \coker (B-I \ C)$ where $(B-I \ C) :
\prod_I \Z \oplus \prod_J \Z \rightarrow \prod_I \Z$.
\end{theorem}

\begin{proof}  Let $F$ be a desingularization of
$E$.  Since $F$ is row-finite and has no sinks it follows
from \cite[Theorem 3.2]{rsz} that $K_0(C^*(E)) \cong \coker
(A_F^t-I)$ and $K_1(C^*(E)) \cong \ker
(A_F^t-I)$. By \cite[Theorem 2.11]{dt-desing} $C^*(E)$ is
Morita equivalent to $C^*(F)$.  Because $K$-theory is
stable, we have that $K_0(C^*(E)) \cong \coker
(A_F^t-I)$ and $K_1(C^*(E)) \cong \ker(A_F^t-I)$.
The result then follows from Lemma
\ref{kerofdesing}.

Furthermore, if $E$ satisfies Condition (L), then it
follows from \cite[Lemma 2.7]{dt-desing} that $F$ also
satisfies Condition (L).  Hence by \cite[Theorem 6.16]{tomf-thesis} we have that $\ext(C^*(F)) \cong \coker (B-I \
C)$.
Since $\ext$ is stable, the result again follows from
Lemma \ref{kerofdesing}.
\end{proof}

\begin{corollary}
If every vertex of $E$ is either a sink or emits infinitely many
edges, then $K_0(C^*(E)) \cong \bigoplus_{E^0} \Z$ and $K_1(C^*(E))
\cong \ext(C^*(E)) \cong \{0\}$.
\end{corollary}

\begin{proof}
$I=\emptyset$, so we have $\bigoplus_I \Z = \prod_I \Z
= \{0\}$, and the result then follows from Theorem \ref{thm-main}.
\end{proof}

In \cite{rsz}, Raeburn and Szyma\'nski prove that every graph
algebra is an Exel-Laca algebra, but not conversely.
In
particular, they produce a matrix
$$
A = \left(\begin{smallmatrix}
1 & 0 & 1 & 1 & 1 & 1 &\\
0 & 1 & 1 & 1 & 1 & 1 &\\
1 & 0 & 1 & 0 & 0 & 0 & \cdots \\ 
0 & 1 & 0 & 1 & 0 & 0 & \\ 
0 & 0 & 1 & 0 & 1 & 0 & \\ 
0 & 0 & 0 & 1 & 0 & 1 & \\
  &   &   & \vdots & & & \ddots \\ 
\end{smallmatrix}\right)
$$
such that the Exel-Laca algebra $\mathcal{O}_A$ is not a graph
algebra.  They do prove, however, that $C^*(E_A)$ is a
$C^*$-subalgebra in $\mathcal{O}_A$, where $E_A$ is the graph
whose vertex matrix is $A$ \cite[Proposition 5.1]{rsz}, and this
prompts them to ask if anything more can be said about the
relationship between the two.  

It appears not.  For 
if $A$ and $E_A$ are as above, the reader can check using Theorem
\ref{thm-main} that $K_0(C^*(E_A)) \cong K_1(C^*(E_A)) \cong \{0\}$.  In
\cite[Remark 4.3]{rsz}, the $K$-theory of $\mathcal{O}_A$ is
computed as $K_0(\mathcal{O}_A) \cong \{0\}$ and
$K_1(\mathcal{O}_A) \cong \Z$.  Hence $C^*(E_A)$ is not a full corner 
of $\mathcal{O}_A$, and in fact $C^*(E_A)$ and $\mathcal{O}_A$ are
not even Morita equivalent. 

We also point out that, for the matrix $A$ above, knowing the
$K$-theory of $C^*(E_A)$
allows one to actually determine
$C^*(E_A)$ up to isomorphism.  
$C^*(E_A)$ is a purely
infinite, simple, separable, nuclear $C^*$-algebra without
unit and hence the Kirchberg-Phillips Classification
Theorem tells us that it it is determined up to Morita
equivalence by its
$K$-theory \cite[Theorem 4.2.4]{kirchphil}.  
Since $\mathcal{O}_2$ has the same
$K$-theory we may conclude that $C^*(E_A)$ is Morita equivalent
to $\mathcal{O}_2$.  Finally, since $E_A$ is
transitive with infinitely many vertices it follows from
\cite[Theorem 2.13]{hjelmborg} that $C^*(E_A)$ is stable.  Hence
$C^*(E_A) \cong
\mathcal{O}_2 \otimes \K$.

\providecommand{\bysame}{\leavevmode\hbox to3em{\hrulefill}\thinspace}


\begin{thebibliography}{10}

\bibitem{bprs}
T.~Bates, D.~Pask, I.~Raeburn, and W.~Szymanski, \emph{The ${C}^*$-algebras of
  row-finite graphs}, New York J. Math. \textbf{6} (2000), 307--324.

\bibitem{CK1}
J.~Cuntz and W.~Krieger, \emph{A class of ${C}^*$-algebras and topological
  {M}arkov chains}, Invent. Math. \textbf{56} (1980), 251--268.

\bibitem{dt-desing}
D.~Drinen and M.~Tomforde, \emph{The ${C}^*$-algebras of arbitrary graphs},
  preprint (2000).

\bibitem{el}
R.~Exel and M.~Laca, \emph{{C}untz-{K}rieger algebras for infinite matrices},
  J. Reine Angew. Math. \textbf{512} (1999), 119--172.

\bibitem{flr}
N.~Fowler, M.~Laca, and I.~Raeburn, \emph{The ${C}^*$-algebras of infinite
  graphs}, Proc. Amer. Math. Soc. \textbf{8} (2000), 2319--2327.

\bibitem{hjelmborg}
J.~Hjelmborg, \emph{Purely infinite and stable ${C}^*$-algebras of graphs and
  dynamical systems}, Ergod. Th. \& Dyn. Sys. (to appear).

\bibitem{kpr}
A.~Kumjian, D.~Pask, and I.~Raeburn, \emph{{C}untz-{K}rieger algebras of
  directed graphs}, Pacific J. Math \textbf{184} (1998), 161--174.

\bibitem{kprr}
A.~Kumjian, D.~Pask, I.~Raeburn, and J.~Renault, \emph{Graphs, groupoids, and
  {C}untz-{K}rieger algebras}, J. Funct. Anal. \textbf{144} (1997), 505--541.

\bibitem{kirchphil}
C.~Phillips, \emph{A classification theorem for nuclear purely infinite simple
  ${C}^*$-algebras}, Doc. Math. \textbf{5} (2000), 49--114.

\bibitem{rsz}
I.~Raeburn and W.~Szymanski, \emph{{C}untz-{K}rieger algebras of infinite
  graphs and matrices}, preprint (1999).

\bibitem{sz}
W.~Szymanski, \emph{Simplicity of {C}untz-{K}rieger algebras of infinite
  matrices}, Pacific J. Math (to appear).

\bibitem{szsemiproj}
\bysame, \emph{On semiprojectivity of ${C}^*$-algebras of directed graphs},
  preprint (2000).

\bibitem{tomf-thesis}
M.~Tomforde, \emph{Computing {E}xt for graph algebras}, preprint (2001).

\end{thebibliography}
\end{document}